\documentclass[10pt]{amsart}
\usepackage{latexsym, amsmath,amssymb}

\renewcommand{\epsilon}{\varepsilon}
\newcommand{\newsection}[1]
{\subsection{#1}\setcounter{theorem}{0} \setcounter{equation}{0}
\par\noindent}

\newtheorem{theorem}{Theorem}

\newtheorem{lemma}[theorem]{Lemma}
\newtheorem{corr}[theorem]{Corollary}

\newtheorem{proposition}[theorem]{Proposition}
\newtheorem{deff}[theorem]{Definition}

\newcommand{\bth}{\begin{theorem}}
\newcommand{\ble}{\begin{lemma}}
\newcommand{\bcor}{\begin{corr}}

\newcommand{\bdeff}{\begin{deff}}

\newcommand{\bprop}{\begin{proposition}}
\newcommand{\ele}{\end{lemma}}
\newcommand{\ecor}{\end{corr}}
\newcommand{\edeff}{\end{deff}}

\newcommand{\eprop}{\end{proposition}}

\newcommand{\la}{\lambda}

\renewcommand{\Pi}{\varPi}

\renewcommand{\epsilon}{\varepsilon}

\begin{document}

\title[Lower bounds for interior nodal sets of Steklov eigenfunctions]
{Lower bounds for interior nodal sets \\ of Steklov eigenfunctions}
\thanks{The first two authors were supported in part by the NSF grant DMS-1361476}
\author{Christopher D. Sogge}
\email{sogge@jhu.edu}
\author{Xing Wang}
\email{xwang@math.jhu.edu}
\author{Jiuyi Zhu}
\email{jzhu43@math.jhu.edu}
\address{Department of Mathematics,  Johns Hopkins University,
Baltimore, MD 21218}

\begin{abstract}
We study the interior nodal sets, $Z_\la$ of Steklov eigenfunctions in an $n$-dimensional relatively
compact manifolds $M$ with boundary and show that one has the lower bounds $|Z_\la|\ge c\la^{\frac{2-n}2}$ for the size of its $(n-1)$-dimensional Hausdorff measure.  The proof is based on a Dong-type identity and estimates for the gradient of Steklov
eigenfunctions, similar to those in \cite{SZ1} and \cite{SZ2}, respectively.
\end{abstract}

\maketitle

\newsection{Introduction}

This article is concerned with lower bounds for the size of nodal sets,
\begin{equation}\label{nodalsets}Z_\la =\{x\in M: e_\la(x)=0 \}, \end{equation}
of real Steklov eigenfunctions in a smooth relatively compact manifold
$(M,g)$ of dimension $n\ge2$ with boundary $\partial M$.  These
eigenfunctions are solutions of the equation
\begin{equation}
	\label{i.1}
	\begin{cases}
	\Delta_g e_\la =0, \quad \text{in } M
	\\
	\partial_\nu e_\la =\la e_\la, \quad \text{on }\partial M,
	\end{cases}
\end{equation}
where $\nu$ is the unit outward normal on $\partial M$.

The Steklov eigenfunctions were introduced by Steklov~\cite{St} in 1902.
They describe the vibration of a free membrane with uniformly distributed
mass on the boundary.  The equation \eqref{i.1} was studied by Calder\'on~\cite{C} as its solutions can be regarded as eigenfunctions of the Dirichlet to Neumann map.

More specifically, the $e_\la$ in \eqref{i.1} satisfy the eigenvalue
problem
$$Pe_\la = \la e_\la,$$
if the Dirichlet to Neumann operator $P$ is defined as 
$$Pf =\partial_\nu Hf|_{\partial M},$$
where for $f\in C^\infty(\partial M)$, $Hf=u$ is the harmonic extension
of $f$ into $M$, i.e., the solution of
\begin{equation*}
\begin{cases}
\Delta_g u(x)=0, \quad x\in M \\
u(x) = f(x), \quad x\in \partial M.
\end{cases}
\end{equation*}

It is well known that $P$ is a self-adjoint classical pseudodifferential operator of order one whose principal symbol agrees with that of the square root of minus the boundary Laplacian on $\partial M$ coming from the metric.  Furthermore, there is an orthonormal basis of real eigenfunctions
$\{e_{\la_j}\}$ such that
$$Pe_{\la_j} = \la_j e_{\la_j}, \quad \text{and }  \quad \int_{\partial M} e_{\la_j}
e_{\la_k} dV_g = \delta_j^k.$$
The spectrum $\la_j$ is discrete, with
$$0=\la_0<\la_1\le \la_2\le \dots, \quad \text{and } \, 
\la_j \to \infty.$$

Recently there has been much work on the study of nodal sets of Steklov
eigenfunctions.  It has largely been focused on the size of the
nodal set 
$${\mathcal N}_\la = \{x\in \partial M: \, e_\la(x)=0\}$$
on the boundary $\partial M$ of $M$.  Bellova and Lin~\cite{BL} proved
that $|{\mathcal N}_\la|\le C\la^6$, if $|{\mathcal N}_\la|$ denotes
 $d-1$ dimensional Hausdorff measure with here $d=n-1$ denoting the dimension of $\partial M$.  Later, Zelditch~\cite{Z} improved these
 results and gave the optimal upper bound $|{\mathcal N}_\la| \le C\la$
 for analytic manifolds using microlocal analysis.
In the smooth case, the last two authors showed in \cite{WZ} showed
that
\begin{equation}\label{aa}
|{\mathcal N}_\la|\ge c\la^{\frac{3-d}2},
\end{equation}
assuming that $0$ is a regular value for $e_\la$.  
This agrees with the best known general lower bounds for
the boundaryless case (see below), but in both \cite{Z} 
and \cite{WZ}  the nonlocal nature of the operators defining
the eigenfunctions presented an obstacle which had to be overcome.

By the maximum
principle, we know that the nodal sets in $M$ must always intersect
the boundary $\partial M$.  In other words, there can be no component
of the nodal set which is closed in $M$.  Thus, it is natural to study
the size of the nodal set in the interior, $M$.  This question was also raised by Girouard and Polterovich in \cite{GP}.

Let us briefly review the literature concerning the study of nodal
sets for compact boundaryless Riemannian manifolds.  Let $\psi_\la$
denote an $L^2$-normalized eigenfunction on of the Laplace-Beltrami
operator on such a smooth $n$-dimensional manifold, i.e.,
$$-\Delta_g \psi_\la = \la^2 \psi_\la. $$
Yau conjectured in \cite{Yau} that one should have
$$c\la \le |Z_\la|\le C\la,$$
if $Z_\la$ denotes the nodal set of $\psi_\la$, and $|Z_\la|$ its
$(n-1)$-dimensional Hausdorff measure.  In the real analytic case
both the upper and lower bounds were established by Donnelly and
Fefferman~\cite{DF}.  The lower bound was established in the
two-dimensional case by Br\"uning~\cite{Br} and Yau (unpublished);
however, in all other cases, the conjecture remains open in the
smooth case.  Recently there has been much work on establishing 
lower bounds in the smooth case when $n\ge 3$.  Colding and Minicozzi~\cite{CM} and then later the first author and Zelditch \cite{SZ1}, \cite{SZ2} showed that
\begin{equation}\label{bb}|Z_\la|\ge c\la^{\frac{3-n}2},
\end{equation}
which matches up with the lower bounds in \eqref{aa} which
were obtained later.  Another proof of \eqref{bb} was given by
Hezari and the second author in \cite{HS}.  

The arguments in 
\cite{SZ1}, \cite{SZ2} and \cite{HS} involved establishing a
Dong-type identity, similar to the one in \cite{Dong}, and then
using either lower bounds for the $L^1$-norms of $\psi_\la$
or upperbounds for its gradient.  We shall use similar arguments
to establish our main result concerning lower bounds for the
$(n-1)$-dimensional Hausdorff measure of the interior nodal sets
of Steklov eigenfunctions contained in the following result.

\begin{theorem}\label{theorem}
Let $M$ be a smooth relatively compact $n$-dimensional manifold
with smooth boundary $\partial M$.  Then there is a constant $c>0$
so that
\begin{equation}\label{steklov}
|Z_\la|\ge c\la^{\frac{2-n}2}
\end{equation}
for the $(n-1)$-dimensional Hausdorff measure of the
nodal sets given by \eqref{nodalsets} of the
Steklov eigenfunctions \eqref{i.1}.
\end{theorem}

We note that this lower bound is off by a half-power versus the best
known lower bounds, \eqref{bb}, for the boundaryless case. 
We shall explain what accounts for this difference after
we complete the proof of Theorem~\ref{theorem}.
 Also,
it seems clear that in the two-dimensional case the lower bound
\eqref{steklov} is far from optimal since the arguments of 
Br\"uning~\cite{Br} and Yau (see also \cite{Savo})  seem to give the optimal
lower bound $|Z_\la|\ge c\la$ using the fact that the nodal
set must intersect any $C\la^{-1}$ ball in $M$
if $C$ is large enough (see e.g. \cite{GP}).

\newsection{An interior Dong-type identity for Steklov eigenfunctions}

As in \cite{SZ1} we shall want to use the Gauss-Green formula
to establish a Dong-type identity which we can use to prove our
lower bound \eqref{steklov}.  We shall be able to do this since
the singular set
$$S_\la =\{x\in \overline{M}: \, e_\la(x)=0 \, \,
\text{and } \, \, \nabla e_\la(x)=0\}$$
is of Hausdorff codimension 2 or more, i.e., $\text{dim }S_\la \le n-2$.
This is true for $S_\la \cap M$ since $e_\la$ is harmonic
in $M$ (see e.g. \cite[Chapter 4]{HL}), while one can, for instance
see that the same is true for $S_\la \cap \partial M$ using the doubling lemma in \cite{Zh}.  In addition, for each $\la$, there are only finitely
many nodal domains (see e.g. \cite{GP}).  Consequently, we may
write $\overline{M}$ as the (essentially) disjoint union
\begin{equation}\label{2.1}
\overline{M}= \bigcup_{i=1}^{k_\la}(D_i^+ \cup Z_i^+\cup Y^+_i)
\, \cup \, \bigcup_{j=1}^{m_\la}(D_j^-\cup Z_j^-\cup Y^-_j),
\end{equation}
where $D^+_i$ and $D^-_j$ are the connected components of 
$\{x\in M: \, e_\la(x)>0\}$ and $\{x\in M: \, e_\la(x)<0\}$,
respectively, while $Z^\pm_k =\partial D^\pm_k\cap M$
and $Y^\pm_k=\overline{D^\pm_k} \cap \partial M$.
Thus,
$$ Z_\la = \bigcup_{i=1}^{k_\la} Z^+_i \, \cup \,
\bigcup_{j=1}^{m_\la}Z^-_j,$$
and
$$\partial M= \bigcup_{i=1}^{k_\la} Y^+_i \, \cup \,
\bigcup_{j=1}^{m_\la}Y^-_j.$$

The boundary of $D^\pm_k$ in $\overline{M}$ is 
$Z^\pm_k \cup Y^\pm_k$.  Since $S_\la$ has codimension
2 or more and $\partial M$ is smooth, we may use the Gauss-Green
formula (see e.g. Theorem 1 on p. 209 of \cite{Evans}) for any $f\in C^\infty(\overline{M})$
to get
\begin{align*}
\int_{D^+_k} \Delta_g f e_\la  \, dV
&=\int_{D^+_k} f\Delta_g e_\la \, dV-\int_{\partial D^+_k}f \partial_\nu e_\la \, dS 
+\int_{\partial D^+_k} \partial_\nu f  e_\la \, dS
\\
&=-\la \int_{Y^+_k}f e_\la dS + \int_{Z^+_k}f|\nabla e_\la|\, dS
+\int_{Y^+_k}\partial_\nu f e_\la \, dS.
\end{align*}
Here $\partial_\nu$ denotes the outward Riemann derivative on 
$\partial D^+_k$, and we used the equation \eqref{i.1} to get the
last equality.  Rearranging, we see from above that 
\begin{equation}\label{22.2}
\la \int_{Y^+_k}f e_\la \, dS 
-\int_{Y^+_k}\partial_\nu f e_\la \, dS
+\int_{D^+_k} \Delta_g f \, e_\la dV=
\int_{Z_k^+}f|\nabla e_\la| \, dS.
\end{equation}
Similarly for each negative nodal domain we have
\begin{align*}
\int_{D^-_k} \Delta_g fe_\la \, dV&=\int_{D^-_k}f \Delta_g e_\la \, dV
-\int_{\partial D^-_k} f \partial_\nu e_\la \, dS
+\int_{\partial D^-_k} \partial_\nu f e_\la \, dS
\\
&=-\la \int_{Y^-_k}f e_\la dS - \int_{Z^-_k}f|\nabla e_\la| \, dS
+\int_{Z^-_k}\partial_\nu f e_\la \, dS,
\end{align*}
using in the last step that on each $Z^-_k$, unlike on each $Z^+_k$,
$\partial_\nu e_\la =|\nabla e_\la|$ since $e_\la$ increases as it
crosses $Z^-_k$ from $D^-_k$.  Rearranging this time leads to
\begin{equation}\label{22.3}
\la \int_{Y^-_k}f e_\la \, dS
-\int_{Y^-_k}\partial_\nu f e_\la \, dS
+\int_{D_k^-}\Delta_g f \, e_\la \, dV
=-\int_{Z^-_k}f|\nabla e_\la|\, dS.
\end{equation}
Since $e_\la >0$ in $D^+_k$ and $e_\la <0$ in $D^-_k$, we can combine
\eqref{22.2} and \eqref{22.3} into
\begin{equation}\label{22.4}
\la \int_{Y^\pm_k}f\, |e_\la| \, dS
-\int_{Y^\pm_k}\partial_\nu f |e_\la| \, dS
+\int_{D^\pm_k}\Delta_g f\, |e_\la|\, dV
=\int_{Z^\pm_k}f|\nabla e_\la|\, dS.
\end{equation}

Since almost every point in $Z_\la$ belongs to exaclty one $Z^+_i$
and one $Z^-_j$ and almost every point in $\partial M$ belongs
to just one of the sets $Y^\pm_k$, if we sum up the identity
\eqref{22.4}, we conclude that we have the Dong-type identity
\begin{equation}\label{22.5}
\la \int_{\partial M}f \, |e_\la|\, dS
-\int_{\partial M}\partial_\nu f \, |e_\la| \, dS
+
\int_M \Delta_g f \, |e_\la|\, dV = 2\int_{Z_\la}f\, |\nabla e_\la| \, dS.
\end{equation}
Of course if $f\equiv 1$ this simplifies to
\begin{equation}\label{22.6}
\la \int_{\partial M} \, |e_\la|\, dS
 = 2\int_{Z_\la} |\nabla e_\la| \, dS,
\end{equation}
which is what we shall use in our proof of Theorem~\ref{theorem}.

\newsection{Interior estimates for Steklov eigenfunctions}

We shall prove interior estimates for the $e_\la$ which are natural
analogs of the ones obtained earlier in the boundaryless case
by Sogge and Zelditch \cite{SZ1}, \cite{SZ2}.  We shall use
arguments which are similar to those of Shi and Xu~\cite{SX} and
\cite{Xu} and H. Smith (unpublished).

Specifically,
we have the following:

\begin{proposition}\label{main}
If $e_\la$ is as above and if $d=d(x)$ denotes the distance from $x\in M$ to $\partial M$,
\begin{equation}
	\label{i.2}
	\|(\la^{-1}+d) \, \nabla_g e_\la \|_{L^\infty(M)}+\|e_\la\|_{L^\infty(M)}
	\le C\la^{\frac{n-2}2}\|e_\la\|_{L^1(\partial M)}.
\end{equation}
\end{proposition}

Let us first argue that on the boundary, we have these estimates. Indeed,
\begin{equation}
	\label{i.3}
	\la^{-\alpha}\|D^\alpha e_\la \|_{L^\infty(\partial M)}
	\le C_\alpha\la^{\frac{n-2}2}\|e_\la \|_{L^1(\partial M)},
\end{equation}
with $D^\alpha$ here referring to $\alpha$ boundary derivatives.  This inequality follows from arguments in \cite{SS} and \cite{SZ1}--\cite{SZ2}, since $Pe_\la =\la e_\la$ where $P$ is a classical self-adjoint pseudodifferential of order one operator whose principal symbol agrees with that of the square root of minus the boundary Laplacian.  As a result we
can use Lemma 5.1.3 in \cite{Sobook} to write $e_\la =T_\la e_\la$, where
$T_\la$ is an integral operator on the $(n-1)$-dimensional boundary of $M$ whose kernel $K_\la(x,y)$ satisfies $D^\alpha K=O(\la^{\alpha +\frac{n-2}2})$ for
each $\alpha$, which immediately gives us \eqref{i.3}.

For the next step, we use that by the maximum principle, the bounds
in \eqref{i.3} for $e_\la$ yield
\begin{equation}
	\label{i.4}
	\|e_\la\|_{L^\infty(M)}\le C\la^{\frac{n-2}2}\|e_\la\|_{L^1(\partial M)},
\end{equation}
as desired.  Thus, we only need to prove the bounds in \eqref{main} for
$\nabla_g e_\la$.

As a first step we realize that we can obtain this estimate 
in the region of $M$ which is of distance $\delta\la^{-1}$ from the boundary
just by using standard Schauder estimates for a given $\delta>0$.  Indeed,
since $e_\la$ is harmonic in $M$ and \eqref{i.3} is valid, it follows from Corollary 6.3 in \cite{GT} 
applied to balls centered at points $x\in M$ or radius $r\le d(x)/2$
that we have
\begin{equation}
	\label{i.5}
	\|d \, \nabla_g e_\la\|_{L^\infty(\{x\in M: \text{dist}(x,\partial M)\ge \delta \la^{-1})}\le C_\delta\la^{\frac{n-2}2}\|e_\la\|_{L^1(\partial M)}.
\end{equation}
Here, the constant $C_\delta$ depends on $\delta$ and $(M,g)$, but not on 
$\la$.

To finish the proof of \eqref{main}, it suffices to show that if $\delta>0$ is sufficiently small we also have the uniform bounds
\begin{equation}
	\label{i.6}
	\la^{-1}\| \nabla_g e_\la\|_{L^\infty(M \cap B(x_0,\delta \la^{-1}))}\le C_\delta\la^{\frac{n-2}2}\|e_\la\|_{L^1(\partial M)}, \quad
	x_0\in \partial M,
\end{equation}
with $B(x_0,\delta \la^{-1})$ denoting the geodesic ball of radius
$\delta \la^{-1}$ about the boundary point $x_0$.

To prove this we shall use local coordinates and a scaling argument.   We shall work in such coordinates and scale and normalize $e_\la$ by replacing
it by 
\begin{equation}\label{i.7}
u_\la(x)=\la^{-\frac{n-2}2}e_\la(x/\la).
	\end{equation}
  Similarly, we shall scale
the $\delta\la^{-1}$ ball so that it becomes a $\delta$ ball 
$\tilde B(x_0,\delta)$ and use the ``stretched" Laplacian
with principal part $\sum g^{jk}(x/\la) \partial_j\partial_k$
(coming from the ``stretched" metric $g_{jk}(x/\la)$), which
denote by $L$.  It follows from \eqref{i.3} that we have
the uniform bounds
\begin{equation}
	\label{i.8}
	\|D^\alpha u_\la \|_{L^\infty(\partial \tilde M
	)}
	\le C_\alpha\|e_\la \|_{L^1(\partial M)},
\end{equation}
where $\tilde M$ denotes the stretched version of $M$ in our local coordinates.  Additionally, the coefficients of our ``stretched" Laplacian $L$ belong to a bounded subset of $C^\infty$ as $\la\ge 1$
and $x_0\in \partial M$ vary.  Also, because of \eqref{i.7} we can
find a function $\varphi_\la$ in our local coordinate system 
which agrees with $u_\la$ on $\partial \tilde M$ and has bounded
$C^{2,\alpha}(\tilde B(x_0,2\delta)\cap \tilde M)$ norm independent of $\la \ge 1$ and $x_0\in \partial M$
for a given $0<\alpha<1$.
Therefore, if we apply Corollary 8.36 in \cite{GT} to $u=u_\la -\varphi_\la$ 
and $f=-L\varphi_\la$,
we conclude that the $C^{1,\alpha}(\tilde B(x_0,\delta))$ norm $u_\la$ is bounded
uniformly with respect to these parameters if $\alpha$ is fixed.  
Thus, we in particular have the uniform bounds
$$\|Du_\la\|_{L^\infty(\tilde B(x,\delta)\cap \tilde M)}
\le C.$$
If we go back to the original local coordinates and recall
\eqref{i.7}, we obtain \eqref{i.6}, which completes
the proof of Proposition~\ref{main}.

\newsection{Conclusion}
It is now very easy to prove Theorem~\ref{theorem}.  If we use
\eqref{22.6} and \eqref{i.2}, we conclude that
$$\la \|e_\la\|_{L^1(\partial M)}=2\int_{Z_\la}|\nabla e_\la|\, dS
\le C\la^{\frac{n-2}2}\|e_\la\|_{L^1(\partial M)}
\int_{Z_\la}(\la^{-1}+d(x))^{-1}\, dS, $$
where, as before, $d(x)$ denotes the distance from $x\in M$ to $\partial M$.  From this, we deduce that
\begin{equation}\label{4.1}
\la^{2-\frac{n}2}\le C\int_{Z_\la}(\la^{-1}+d(x))^{-1} \, dS.
\end{equation}
Clearly this inequality yields \eqref{steklov}, establishing Theorem~\ref{theorem}.

\bigskip
\noindent
{\bf Remarks:}  There is a simple explanation of why the lower bounds
\eqref{steklov} are off by a half power versus the corresponding best
lower bounds \eqref{bb} for the boundaryless case.  This is because
the Dong-type identity in \cite{SZ1} involved $\la^2$ in the left
side instead of $\la$, which accounts for a relative loss of a full power
of $\la$, but, on the other hand, the estimates for the gradient
here are one half power better due to the fact that the boundary
of $M$ is of one less dimension, accounting for a relative gain
of a half power.

In some cases one can use \eqref{4.1} to get improved lower bounds.  For
instance if we let
$$Z_{\la,k}=\{x\in Z_\la: \, d(x)\in [2^{-k}, 2^{-k+1})\}$$
and if $|Z_{\la,k}|\le C2^{-k}|Z_\la|$ for $C\le k\le \log_2 \la$
and if $|\{x\in Z_\la: d(x)\le \la^{-1}\}|\le C\la^{-1}|Z_\la|$, with
$C$ fixed, we then
get the lower bound $|Z_\la|\ge c \la^{2-\frac{n}2}/\log\la$, which is
essentially optimal when $n=2$.  The subsets $Z_{\la,k}$ of $Z_\la$ have this property, for instance, for the Steklov eigenfunctions
$r^m \sin m\theta$ on the disk in ${\mathbb R}^2$ (written in polar coordinates).


\begin{thebibliography}{MA}
\bibitem{BL} K. Bellova and F.H. Lin:{\em  Nodal sets of Steklov eigenfunctions},
arXiv:1402.4323.

\bibitem{Br} J. Br\"uning: {\em  \"Uber Knoten von Eigenfunktionen des Laplace-Beltrami-Operators},
(German)  Math. Z.  {\bf 158}  (1978),  15--21.

\bibitem{C} A.P. Calder\'{o}n: {\em  On a inverse boundary value problem}, in ``Seminar
in Numerical Analysis and Its Applications to Continuum Physics",
1980, 65--73, Soc. Brasileira de Matem\'{a}tica, Rio de Janeiro.

\bibitem{CM} T.H. Colding and W. P. Minicozzi II:{\em  Lower bounds for nodal sets
of eigenfunctions}, Comm. Math. Phys. {\bf 306} (2011), 777--784.

\bibitem{Dong} R-T Dong: {\em Nodal sets of eigenfunctions on Riemann surfaces},
 J. Differential Geom.  {\bf 36}  (1992),  493--506.

\bibitem{DF} H. Donnelly and C. Fefferman: {\em Nodal sets of eigenfunctions
on Riemannian manifolds}, Invent. Math. {\bf 93} (1988), 161--183.

\bibitem{Evans} L.C. Evans and R. Gariepy: {\em Measure theory and fine properties of functions}, CRC Press, Boca Raton, Ann Arbor, and London, 1992.

\bibitem{GT} D. Gilbarg and N. S. Trudinger:{\em   Elliptic partial differential equations of second order},
Reprint of the 1998 edition.
Classics in Mathematics. Springer-Verlag, Berlin,  2001.

\bibitem{GP} A. Girouard and I. Polterovich: {\em Spectral geometry of the Steklov
problem}, arXiv:1411.6567.

\bibitem{HL} Q. Han and F.H. Lin: {\em Nodal sets of solutions of elliptic
differential equations}, book in preparation (online at
http://www.nd.edu/qhan/nodal.pdf).

\bibitem{HS} H. Hezari and C.D. Sogge: {\em A natural lower bound for the
size of nodal sets}, Anal. PDE. {\bf 5} (2012)  1133--1137.


\bibitem{Savo} A. Savo: {\em Lower bounds for the nodal length of eigenfunctions of the
 Laplacian},
 Ann. Global Anal. Geom.  {\bf 19}  (2001),  133--151.

\bibitem{SS} A. Seeger and C.D. Sogge: {\em Bounds for eigenfunctions
of differential operators}, Indiana Math. J. {\bf 38} (1989), 669--682.

 \bibitem{SX} Y. Shi,  and Xu, Bin: {\em  Gradient estimate of an eigenfunction on a compact Riemannian manifold
 without boundary},
 Ann. Global Anal. Geom.  {\bf 38}  (2010),  21--26.

\bibitem{Sobook} C.D. Sogge: {\em Fourier integrals in classical analysis},
Cambridge Tracts in Mathematics, {\bf 105} Cambridge University Press,
Cambridge, 1993.

\bibitem{S1} C.D. Sogge: {\em  Concerning the $L^p$ norm of spectral clusters for
second-order elliptic operators on compact manifolds}, J. Funct.
Anal. {\bf 77} (1988), 123--138.

\bibitem{St} W. Stekloff: {\em  Sur les probl\`{e}mes fondamenbetax de la
physique mathematique}, \emph{Ann. Sci. \'{E}cole Norm. Sup.}
{\bf 19} (1902), 191-259.

\bibitem{SZ1} C.D. Sogge and S. Zelditch: {\em Lower bounds on the
Hausdorff measure of nodal sets}, Math. Res. Lett. {\bf 18} (2011), 25--37.

\bibitem{SZ2}C.D. Sogge and S. Zelditch: {\em   Lower bounds on the Hausdorff measure of nodal sets II},
Math. Res. Lett. {\bf 19}  (2012),  1361--1364.

\bibitem{WZ} X. Wang and J. Zhu: {\em  A lower bound for the nodal sets of Steklov
eigenfunctions}, arXiv:1411.0708.

 \bibitem{Xu} X. Xu: {\em  Gradient estimates for the eigenfunctions on compact manifolds with
 boundary and H\"ormander multiplier theorem},
 Forum Math.  {\bf 21}  (2009), 455--476.
 
 \bibitem{Yau} S. T. Yau: {\em Survey on partial differential equations in differential geometry}, Seminar on Differential Geometry, pp. 3–71, Ann. of Math. Stud., {\bf 102}, Princeton Univ. Press, Princeton, N.J., 1982.

\bibitem{Z} S. Zelditch: {\em  Measure of nodal sets of analytic
Steklov eigenfunctions}, arXiv:1403.0647.

\bibitem{Zh} J. Zhu: {\em Doubling property and vanishing order of Steklov
eigenfunctions},   Comm. Partial Differential Equations.
\end{thebibliography}
\end{document}